\documentclass[12pt]{amsart}
\usepackage{latexsym}
\theoremstyle{plain}

\theoremstyle{definition}

\theoremstyle{remark}

\def\R{{\rm I\kern-2ptR}}
\def\N{{\rm I\kern-2ptN}}

\def\wider{\rhd}
\def\narrower{\lhd}

\def\qed{\hskip .6em \raise1.8pt\hbox{\vrule height4pt width6pt depth2pt}}

\begin{document}
\baselineskip 17.9pt

\centerline {\bf A characterization of operators preserving disjointness
in terms of their inverse}

\bigskip

\centerline {\bf Y. Abramovich and A. Kitover}

\bigskip
\hskip 1.3cm {\bf Abstract.}
The characterization mentioned in the title is found.

\bigskip
{\bf 1. Introduction.}
Recall that a (linear) operator $T:X\to Y$ between vector
lattices is {\bf disjointness preserving} if $T$ sends elements disjoint in $X$
to elements disjoint in $Y$. If $T$ is a bijective disjointness preserving 
operator between Banach lattices, then a well known theorem by 
Huijsmans--de Pagter [7] and Koldunov [8]  asserts that the inverse
$T^{-1}: Y\to X$ is also  disjointness preserving. 
Many other results describing various conditions under which $T^{-1}$
is disjointness preserving
can be  found in [5]. It was believed  for a while
that the same conclusion  should  remain true for disjointness    
preserving    operators between arbitrary vector lattices, or at least, for
operators between Dedekind complete vector lattices. 
However,
as has been recently shown  by the authors [4,5], this is      
not true  in general. This means, in particular, that if one wants to       
find a characterization of a disjointness preserving  operator
in terms of its inverse, then a different condition is needed
rather than disjointness preservation. 

It is the purpose of this note  to present such a condition.                 
The authors would like to express their  thanks to Beata       
Randrianantoanina for her help in identifying this condition.
In her talk\footnote{Delivered at the conference
          {\it Function Spaces}, held in Edwardsville in May of
1998.} devoted to description of  non-surjective     
isometries between  some Orlicz  spaces and based on her work [9], 
Randrianantoanina  introduced an interesting monotonicity  condition         
and asked if it implied disjointness preservation.
The essence of this condition is as follows: 
if the support of a measurable function $x_1$
is contained in the support of another measurable function $x_2$, 
then the same is true
for the supports of their images, that is, 
the support of  $Tx_1$ is contained in the
the support of $Tx_2$, where $T$ is the isometry in question.
An abstract order-theoretic version 
of this condition will be introduced in Definition~2.2
and denoted by $(\beta)$. As  Examples~2.5
and~2.6 demonstrate, condition
$(\beta)$  and disjointness preservation are independent in general.
Nevertheless, these conditions are related in more than one way.
First of all, as  Theorem~3.4 demonstrates, it is precisely   
condition  $(\beta)$ that characterizes the inverses of
operators preserving  disjointness. 
This characterization allows us  to describe 
bijective operators for which         
condition $(\beta)$  and disjointness preservation are
equivalent. Roughly speaking, they are  equivalent
if and only if $T^{-1}$ is
disjointness preserving. The most important instances of this
situation
are presented in Theorem~3.7. Accordingly, each example 
of a bijective disjointness preserving operator $T$ whose inverse       
$T^{-1}$ is not  disjointness preserving, is an example of a 
disjointness preserving operator that does not satisfy condition      
$(\beta)$.

In our terminology regarding vector lattices
and operators on them we follow [6]. Throughout the work all vector 
lattices are assumed to be Archimedean. The reader is referred to [5]
for a comprehensive study of  the inverses of operators preserving   
disjointness.

\bigskip

{\bf 2. A sufficient condition.}
 Recall that for a subset $A$ of a vector lattice $X$
the symbol $A^d$ denotes the disjoint complement of $A$ which is
defined as follows:
$A^d=\{x\in X: x\perp a \ \ {\rm for \ each}\ \ a\in A\}$.
The disjoint complement $(A^d)^d$ of $A^d$ is denoted simply by $A^{dd}$;
the set $A^{dd}$ is the band generated by $A$.
\bigskip

\noindent
{\bf Definition 2.1.}
{\sl Let $a,b$ be two elements in a vector lattice $X$. 
We  say that $a$ and $b$ are  {\bf of the same width} 
if $\{a\}^{dd} = \{b\}^{dd}$, that is, $a$ and $b$
generate the same band in $X$. 
Similarly, we say that $b$ is {\bf wider} than  $a$ if 
$\{b\}^{dd} \supseteq \{a\}^{dd}$. 
 }
\medskip

Clearly, $b$ is {\bf wider} than  $a$ if and only if $\{b\}^{dd} \ni a$.
If $b$ is wider than $a$, we will write $b  \wider a$. 
An equivalent notation is $ a \narrower b$, in which case
we will also say that $a$ is narrower than $b$.

\bigskip
\noindent
{\bf Definition 2.2.}
{\sl Let us say that a $($linear$)$ operator $T:X\to Y$ between           
vector lattices satisfies condition $(\beta)$ if $Ta \narrower Tb$ in $Y$
whenever $a \narrower b$ in $X$.
}

Note that  if $T:X\to Y$ satisfies condition $(\beta)$, then
{\sl for each $x\in X$ the operator $T$ sends the band $\{x\}^{dd}$
into the band $\{Tx\}^{dd}$.}

Indeed,
take an arbitrary element $u\in \{x\}^{dd}$ and show that 
$Tu \in \{Tx\}^{dd}$.
Note  that the former condition $u\in \{x\}^{dd}$ is equivalent to
saying that $u\narrower x$ and hence, in view of $(\beta)$, we have
$Tu\narrower Tx$. This means that  $Tu\in \{Tx\}^{dd}$.

In particular,  if $T$  satisfies  $(\beta)$ and $Tx=0$ for some $x\in X$,
then $Tx'=0$ for each element $x'$ in the band generated by $x$. That is, 
together with each $x$ the kernel, $ker (T)$, of the operator $T$ contains the
principal band generated by $x$, and thus $ker (T)$ is an ideal.

\bigskip

\noindent
{\bf Corollary 2.3.}
{\sl Let $T:X\to Y$ satisfy condition $(\beta)$.
For each $x\in X$ we have $T(|x|) \in \{Tx\}^{dd}$.}
\medskip

\bigskip
For any  injective operator $T: X\to Y$ between vector lattices we denote 
by $T^{-1}$ the inverse operator  defined on $TX$.

\bigskip
\noindent
{\bf Theorem 2.4.}
{\sl Let $T:X\to Y$ be an injective  operator between 
vector lattices. If $T$ satisfies condition $(\beta)$, then
$T^{-1}: TX\to X$ is a disjointness preserving operator.}

{\bf Proof.}
Take arbitrary disjoint elements  $y_1,y_2$ in $TX$ and let 
$x_i:=T^{-1}y_i, \ i=1,2$. We want to show that $x_1\perp x_2$. Let
$v_i=T(|x_i|), \ i=1,2$. By the previous corollary we know that $v_i\in
\{y_i\}^{dd}$, whence
$v_1 \perp v_2$.

Consider $u = |x_1|\land |x_2|$. Certainly $u\narrower |x_i|$ for $i=1,2$ 
and so,  in view of $(\beta)$, we have $Tu \narrower T(|x_1|)= v_1$. Thus,
$Tu \in \{v_1\}^{dd}$. Similarly, we have 
$Tu \in \{v_2\}^{dd}$. These two inclusions and the fact that 
$v_1\perp v_2$ imply that $Tu=0$, that is, $u=0$ because $T$ 
is injective. We are done since $u=|x_1|\land |x_2|$.    
\qed

\bigskip

It is interesting to point out that Theorem~2.4 does not claim that the  
subspace 
$TX$ of $Y$ necessarily has disjoint elements. 
It only claims that if  they do exist then $T^{-1}$ sends them to disjoint   
elements. Let us consider a simple example of an operator $T$ satisfying
condition
$(\beta)$ and such that 
$TX$ does not have non-trivial disjoint elements.

\bigskip
\noindent
{\bf Example 2.5.} Let $X=\R^2$, $Y=C[0,1]$ and  functions $v_1, v_2 \in Y$ be
defined by
$v_1(t) = 1, v_2(t) =t \ \ {\rm for}\ t\in [0,1]$.
Let $Te_i=v_i$, where $e_1$ and $e_2$ are the standard unit vectors in $X$. 

It is
obvious that for each non-zero $x\in X$ the function $y=Tx$ has full support
in $Y$, and so $T$ satisfies $(\beta)$. At the same time $TX$ does not  
contain any non-trivial pair of disjoint  elements.
\qed

\medskip
Clearly, the operator $T$ constructed in Example~2.5 is not   
disjointness preserving,
and so, in particular, {\it condition $(\beta)$ does not imply
disjointness preservation}. Another kind of examples with the latter
property is  provided by any integral operator with a positive
kernel. To see  that the converse implication does not hold either
(in other words, to see that
{\it $(\beta)$ and disjointness preservation are  independent})
we need to present also an example of {\it a disjointness preserving 
operator which does not satisfy  condition $(\beta)$}. 
As said earlier, the shortest way to accomplish this is to use any            
of the examples constructed in [4,5] of a bijective disjointness      
preserving operator $T$ between vector lattices (which can be, even,    
Dedekind complete) such that $T^{-1}$ is not disjointness preserving.         
In view of Theorem~2.4 such
$T$ cannot satisfy  $(\beta)$. However, none of these counterexamples is too
simple. Therefore, to make our work self-contained, we
will  present below  an independent  example 
demonstrating that a disjointness preserving 
operator  does not imply $(\beta)$ in general.

\bigskip
\noindent
{\bf Example 2.6.} There exists a disjointness preserving
operator $T: X\to Y$ between vector lattices 
such that $T$ fails condition $(\beta)$.

Let
$X_1$ be the collection of all continuous functions
on $[0,2]$ that vanish at some neighborhood of $[1,2]$,
that is, 
$$
X_1=\{x\in C[0,2] : \ \ \exists \varepsilon_x > 0\ \ {\rm such \ that}\ \
x(t) =0 \ {\rm for  \ each}\ \ t > 1-\varepsilon_x. \}
$$
Similarly,
$$
X_2=\{x\in C[0,2] : \ \ \exists \varepsilon_x > 0\ \ {\rm such \ that}\ \
x(t) =0 \ {\rm for  \ each}\ \ t < 1+\varepsilon_x. \}
$$
Let $X$ be the linear span of $X_1$, $X_2$ and the constant function $\bf 1$.
It is easy to verify that $X$ is a vector sublattice of  $C[0,2]$.
For each function $x$ in $X$ we can find
unique $x_i\in X_i$ and a scalar $\lambda \in \R$ such that
$x$ has the following representation:
$x=x_1+x_2 +\lambda \bf 1$.

Let $Y=X_1\oplus X_2\oplus \R$ be the order direct sum  of $X_1$, $X_2$      
and
\R. Thus for each  $y$ in $Y$ we can find
unique $y_i\in X_i$ and a scalar $\lambda \in \R$ such that
$y= (x_1, x_2,\lambda)$.

Now we are ready to define an operator $T: X\to Y$ by letting
$Tx =(x_1, x_2,\lambda)$.
Since $T({\bf 1})= (0, 0, 1)$ it follows obviously 
that $T$ fails condition $(\beta)$. It remains to verify that $T$ is
disjointness preserving. Take any disjoint  $x', x'' \in X$.
Then $x'=x'_1+x'_2 +\lambda' \bf 1$
and $x''=x''_1+x''_2 +\lambda'' \bf 1$. 
A crucial thing to observe  now is the fact that $x'\perp x''$ if and only if
$\lambda' =\lambda''=0 $, $\ x'_1 \perp x_1''$ and 
$x'_2 \perp x_2''$. Thus
$$
Tx' =(x'_1, x'_2 ,0) \ \ {\rm and} \ \ Tx'' =(x''_1, x''_2, 0),
$$
whence $Tx' \perp Tx''$.
\qed

It is worth pointing out that the 
counterexamples from [4,5] mentioned above 
are possible   for non-regular    
operators only. A simple verification shows that the disjointness preserving  
operator $T$ constructed in Example~2.6 and  failing $(\beta)$  is not regular
either. This is not accidental as the following important result
reveals.

\bigskip

\noindent
{\bf Theorem 2.7.}
{\sl Let $T$ be a disjointness preserving  operator
between vector lattices. If $T$ is regular, then $T$
satisfies condition $(\beta)$.}

{\bf Proof.}
We will present  only a sketch of the proof, since it
depends on a theorem (Theorem~3.4 in [3]) which is rather far from       
our discussions here. The essence of that  theorem  is that it
allows one to represent each regular disjointness preserving            
operator as a  weighted composition operator. And, for the latter          
class of  operators, condition $(\beta)$ can be verified directly.
\qed

\bigskip

\noindent
{\bf Corollary 2.8.}
{\sl Let $T$ be a  continuous operator
between normed vector lattices. If $T$ is  disjointness preserving, 
then $T$ satisfies condition $(\beta)$.}

{\bf Proof.} As shown in [1, Corollary~1] each {\it continuous} disjointness 
preserving operator between normed vector lattices is necessarily regular,
and so Theorem~2.7 is applicable.
\qed

\bigskip

{\bf 3. A necessary condition.}
Recall that an element $x'$ of a
vector lattice $X$  is said to be  a {\bf component} of an element $x \in X$  
if
$|x'|\land |x-x'|=0$. Our next definition introduces a very weak property 
describing vector lattices whose elements have  relatively large amounts of  
components.

\bigskip
\noindent
{\bf Definition 3.1.}
{\sl We say that a            
vector lattice $X$ has {\bf sufficiently many components} if
whenever $x\notin \{u\}^{dd}$ for some $x,u\in X$ there exists 
a non-zero component $x'$ of $x$ such that $x'\perp u$. }

\medskip 
It is obvious that each vector lattice with the principal    
projection property, or even with a cofinal family of band-projections [5]
has sufficiently many components. In particular, each Dedekind complete   
vector lattice has sufficiently many components.

\bigskip

\noindent
{\bf Theorem 3.2.}
{\sl Let $T:X\to Y$ be an injective  operator between vector lattices,
and assume that $X$ has  sufficiently many components.
If $T$ is a disjointness preserving operator, then $T^{-1}: TX\to X$
satisfies condition $(\beta)$.}

{\bf Proof.}  Fix any $y_0=Tx_0 $ in $Y$ and  take any $y\in TX$
which is narrower than $y_0$.  So $Tx=y$ for some $x\in X$.
We want to show that $x \narrower x_0$. If  not, then 
using the hypothesis that $X$ has sufficiently many components 
we can find a non-zero component $x'$ of $x$ that is disjoint to
$x_0$.

Since $T$ preserves disjointness,  we know that  $T$ sends components to
components,   and so $Tx'$ is a component of $y=Tx$, in particular, 
$Tx'\in \{Tx\}^{dd}=\{y\}^{dd}$. But $y$ is narrower than $y_0$
and, consequently, $Tx'\in \{Tx_0\}^{dd}=\{y_0\}^{dd}$.     
Note that $Tx' \neq 0$ as $x'\neq 0$ and $T$ is injective.
On the other hand, $x' \perp x_0$ implies that  $Tx' \perp Tx_0$,
a contradiction.
\qed

\bigskip
We do not know whether or not the assumption in Theorem~3.2 that 
$X$ has sufficiently many components is essential. As we show next,
under an additional hypothesis about the operator
we can get rid of that assumption. 
Note, however, that if this assumption is indeed essential,
then   it will be rather difficult to produce
a counterexample to this effect.  
One possible approach to this question may be related to
Problem~P.4.2 in [5].

Recall that a
bijection  $T:X\to Y$ between vector lattices is a {\bf d-isomorphism}
if both $T$ and $T^{-1}$ preserve disjointness.

\bigskip

\noindent
{\bf Proposition 3.3.}
{\sl If $T:X\to Y$ is a d-isomorphism, where $X$ and $Y$ are arbitrary vector
lattices, then
$T^{-1}$ satisfies condition $(\beta)$.}

{\bf Proof.} 
Let $y_1=Tx_1, y_2=Tx_2$ be two elements in $Y$ such that
$y_2$ is wider than $y_1$. Assume, contrary to what we want to prove, that
$x_2$ is not wider than $x_1$, that is, $x_1 \notin \{x_2\}^{dd}$.
Therefore there exists some non-zero $x\in X$ such that $x\perp x_2$
and $x\not\perp  x_1$. Since $T$ preserves disjointness 
and $x\perp x_2$ we have
$Tx\perp Tx_2$, and since  $x\not\perp  x_1$ we have
$Tx\not\perp  Tx_1$ (here we use the fact that each d-isomorphism
sends non-disjoint elements to non-disjoint).
Thus we have that $Tx \perp y_2$ and $Tx\not\perp y_1$. This
contradicts our assumption that
$y_2$ is wider than $y_1$.
\qed

\bigskip
Combining Theorems~2.4 and~3.2 and Proposition~3.3 we obtain
immediately a description of duality between condition $(\beta)$ and
disjointness  preservation. 

\bigskip
\noindent
{\bf Theorem 3.4.}
{\sl Let $T:X\to Y$ be a bijective  operator between vector lattices,
and assume that $X$ has sufficiently many components. Then
$T$ is disjointness preserving  if and only if  $T^{-1}$
satisfies condition $(\beta)$.
}

\bigskip

Since each Banach function space is necessarily Dedekind complete,
it certainly has  sufficiently many components. Therefore, our     
next result is an immediate consequence of the previous theorem.
We have singled  this case out since it may be of special interest in
dealing with the isometric operators on Banach function spaces.

\bigskip
\noindent
{\bf Corollary 3.5.}
{\sl Let $T$ be a bijective operator between Banach function
spaces. Then $T$ is disjointness preserving if and only if 
$T^{-1}$ satisfies condition $(\beta)$.
}

\bigskip

As shown earlier,  condition $(\beta)$ and
disjointness preservation are not equivalent in general.
However, there are many situations when they are,
and, as Theorems~2.4 and~3.4 show,
all  these cases reduce to those which
guarantee that $T^{-1}$ is disjointness preserving
when $T$ is.
Theorem~3.7 singles out some cases that are  most 
important for applications. We precede this theorem with a useful    
proposition showing that for a disjointness preserving $T$, 
the operator $T$ satisfies $(\beta)$  if and only if
$T^{-1}$ is disjointness preserving. 

\bigskip

\noindent
{\bf Proposition 3.6.}
{\sl For a disjointness preserving bijection  
$T$ between vector lattices the following two    
statements are equivalent.

{\rm 1)}  $T$ satisfies condition $(\beta)$.

{\rm 2)}  $T^{-1}$ is disjointness preserving.
}

{\bf Proof.} 
The implication $1) \Longrightarrow 2) $ is valid by Theorem~2.4.
Conversely, assume that 2) holds, and thus
$T$ is a d-isomorphism. It remains to  apply Proposition~3.3 to the
operator $T= (T^{-1})^{-1}$ so that we can conclude that $T$ satisfies
$(\beta)$.
\qed

\bigskip
\noindent
{\bf Theorem 3.7.}
{\sl Let $T: X\to Y$ be a bijective  operator between
vector lattices that satisfy any one of the following
non-exclusive conditions:

{\rm 1)} $X$ and $Y$ are Banach lattices.

{\rm 2)} $X$ is $(r_u)$-complete and $Y$ is a normed vector lattice.

{\rm 3)} $X$ is a normed Dedekind $\sigma$-complete vector lattice  and   $Y$ is
$(r_u)$-complete.

{\rm 4)} $X$ is a Dedekind $\sigma$-complete normed vector lattice
and $Y=X$.

\noindent Then $T$ is disjointness preserving if and only if it satisfies
condition $(\beta)$.
} 

{\bf Proof.} We will consider only case 1), which is the most important.
The other cases can be dealt with  similarly. Let $X$ and $Y$ be Banach lattices
and $T:X\to Y$ be disjointness preserving. Hence by the Huijsmans--de
Pagter--Koldunov Theorem the inverse operator
$T^{-1}:Y \to X$ is also disjointness preserving, that is,
$T$ is a d-isomorphism.
Therefore, by  Proposition~3.3, the operator $T= (T^{-1})^{-1}$
satisfies  condition $(\beta)$.

Conversely, assume that $T$ satisfies
condition $(\beta)$. Then by Theorem~2.4 the inverse $T^{-1}: Y \to X$
is disjointness preserving, and the second application of 
the Huijsmans--de Pagter--Koldunov Theorem guarantees that 
$T= (T^{-1})^{-1}$ is also disjointness preserving.
\qed 

We would like to emphasize that the operators, we were dealing
with in  Secion~3, were not assumed to be continuous. Whenever they are,     
the proofs can be simplified in view of Corollary~2.8.
\bigskip

{\bf 4. Some concluding results.}
We are going to address now a natural
question  on the relationship between property $(\beta)$
and a similar property in which
one considers only the elements of the same width     
instead of elements subjected to  wider/narrower conditions.

\bigskip
\noindent
{\bf Definition 4.1.}
{\sl Let us say that an operator $T:X\to Y$ between           
vector lattices satisfies condition $(\beta_0)$ if $T$ sends                 
any two elements of the same width to  elements of the same width.
}

\medskip
It is obvious that condition $(\beta)$ implies $(\beta_0)$. We show next
that the converse implication is also true.

\bigskip
\noindent
{\bf Theorem 4.2.}
{\sl Conditions  $(\beta)$ and  $(\beta_0)$ are equivalent}

{\bf Proof.}
Take any  $a,b \in X$ such that $a  \wider b$, that is,
$a$ is wider than $b$. Observe that for each $x\in X$
the elements $Tx$ and $T(|x|)$ are, 
in view of $(\beta_0)$,  of the same width.
Therefore, without loss of generality, we may assume that           
both $a$ and $b$ are positive.

Consider the element $b+\frac{1}{n} a$, $n=1,2\ldots$.
Since $a,b\ge 0$ and $a \wider b$, it is obvious that
the elements $a$ and $b+\frac{1}{n} a$ are of the same width.
Hence, by  $(\beta_0)$, the images 
$Ta$ and $T(b+\frac{1}{n}a)= Tb+ \frac{1}{n} Ta$ are 
also of the same width.

Finally note  that  the sequence $\{Tb+ \frac{1}{n} Ta\}_n$ order
converges (in actuality,
$(r_u)$-converges) to $Tb$.  This and the fact that
$Ta$ and $ Tb+ \frac{1}{n} Ta$ are of the same width
imply that $Ta \wider Tb$. 
\qed

\bigskip

We proceed to observe that if we consider only the positive
elements in condition $(\beta)$, then, surprisingly enough, the resulting
condition is not equivalent to $(\beta)$. To make all this precise
we introduce a formal  definition.

\bigskip
\noindent
{\bf Definition 4.3.}
{\sl Let us say that an operator $T:X\to Y$ between           
vector lattices satisfies condition $(\beta_+)$ if $Ta \narrower Tb$ in $Y$
whenever $a \narrower b$ and $a,b\in X_+$.
}
\bigskip

It is obvious that $(\beta) \Rightarrow (\beta_+)$. Our next example shows
that the converse is not true in general.

\bigskip
\noindent
{\bf Example 4.4.}
There exists an operator satisfying $(\beta_+)$ but not
satisfying $(\beta)$.

Let $X =L^\infty[-1,1]$ and let $T: X\to X$ be defined as follows:

$$
Tx(t)=\left\{ 
 \begin{array} {ll}
  0 & \ \ \textrm{if $t\le 0$} \\
  x(t) + x(-t)  & \ \ \textrm{if $t \ge 0$}.\\
 \end{array} \right.
$$
Obviously $T\ge 0$. A
straightforward verification shows that $T$ satisfies
$(\beta_+)$. 

However, $T$ fails $(\beta)$. Indeed, consider $e = -\chi_{[-1,0]} +
\chi_{[0,1]} \in X$. Clearly $Te=0$. 
Consequently, if $T$ satisfied $(\beta)$, then, by the comments preceding
Corollary~2.3, the kernel of $T$ would contain the principal band
generated by $e$.  Since $e$ has 
full support in $X$, the band
generated by $e$ coincides with $X$.
But $T$ is not identically zero, a
contradiction. 
\qed

\bigskip

We conclude with one more remark. In hindsight,
an antecedent of  condition $(\beta)$  can be traced in some         
earlier work. Namely, in [2,3] the authors considered the disjointness
preserving operators sending weak units to weak units. In the terms of the
present work this   can be expressed by saying that the elements with
full    support are being mapped to elements also with full support.          
In other words, operators satisfying $(\beta_0)$ form a special subclass
of operators considered in [2,3].

\bigskip
\centerline {\bf References}
\bigskip

[1] Y. A. Abramovich, Multiplicative representation of operators 
        preserving disjointness, { \it Netherl. Acad. Wetensch. Proc. 
                      Ser. A  \bf 86}$\,$(1983), 265--279.

[2] Y. Abramovich, E.~Arenson  and A.~Kitover,
   Operators in Banach C(K)-modules and their spectral properties, 
   {\it Soviet Math. Dokl.\/ \bf 38}$\,$(1989), 93--97.

[3] Y. Abramovich, E.~Arenson  and A.~Kitover,
   {\it Banach $C(K)$-modules and operators preserving disjointness},     
    Pitman Research Notes in
    Mathematical Series $\#277$, Longman Scientific $\&$ Technical, 1992.

[4] Y. A. Abramovich and A. K. Kitover, A solution to a problem on
         invertible disjointness preserving operators,
         {\it Proc. Amer. Math. Soc.\/ \bf 126}$\,$(1998), 
         1501--1505.

[5] Y. A. Abramovich and A. K. Kitover, {\it Inverses of disjointness 
      preserving operators}, {\it Memoirs of the Amer. Math. Soc.}, 
      forthcoming.

[6] C. D. Aliprantis and O. Burkinshaw, {\it Positive Operators\/}, 
        Academic Press, New York \& London, 1985.

[7] C. B. Huijsmans and  B. de Pagter, Invertible disjointness preserving
       operators, {\it Proceed. Edinburgh. Math. Soc. {\rm (2)} \bf 37}$\,$   
       (1993), 125--132.

[8] A. V. Koldunov, Hammerstein operators  preserving   disjointness,
       {\it Proc. Amer. Math. Soc.\/ \bf 123}$\,$(1995), 1083--1095.

[9] B. Randrianantoanina, Injective isometries in Orlicz spaces, 1998,
preprint.

\bigskip

Y. A.  Abramovich                    \hskip 5.5cm       A. K. Kitover

Department of Mathematical Sciences \hskip 1.94cm    Department of Mathematics

IUPUI, Indianapolis, IN 46202    \hskip 3.3cm     CCP, Philadelphia, PA 19130

USA                               \hskip 8cm    USA

yabramovich@math.iupui.edu   \hskip 3.55cm   akitover@ccp.cc.pa.us

\end{document}